\documentclass[12pt]{amsart}
\usepackage{amsmath}
\newtheorem{thm}{Theorem}[section]
\newtheorem{theorem}[thm]{Theorem}
\newtheorem{corollary}[thm]{Corollary}
\newtheorem{lemma}[thm]{Lemma}
\newtheorem{proposition}[thm]{Proposition}
\newenvironment{theorem*}[1]{\smallskip\noindent{\bf #1}\em}
                        {\medskip\rm}

\theoremstyle{definition}

\theoremstyle{remark}
\newtheorem{remark}[thm]{Remark}

\numberwithin{equation}{section}

\newcommand\ov{\overline}

\newcommand\wt{\widetilde}
\newcommand\wh{\widehat}

\def\oast{\,\ov{\vphantom{g}\ast}\,}
\def\hast{\,\widehat\ast\,}
\def\tast{\,\widetilde\ast\,}

\newcommand\al{\alpha}
\newcommand\be{\beta}
\newcommand\ga{\gamma}
\newcommand\la{\lambda}
\newcommand\si{\sigma}
\renewcommand\th{\theta}
\newcommand\om{\omega}

\newcommand\bC{{\mathbb C}}
\newcommand\bN{{\mathbb N}}
\newcommand\bR{{\mathbb R}}

\newcommand\cW{{\mathcal W}}

\newcommand\fA{{\mathfrak A}}
\newcommand\fD{{\mathfrak D}}
\newcommand\fm{{\mathfrak m}}

\begin{document}

\title[Inverse spectral problems]{Inverse spectral problems
    for Sturm-Liouville operators with singular potentials, IV.\\
    Potentials in the Sobolev space scale${}^{\dag}$}%

\author[R.~O.~Hryniv, Ya.~V.~Mykytyuk]{Rostyslav O.~Hryniv and Yaroslav V.~Mykytyuk}%
\address{Institute for Applied Problems of Mechanics and Mathematics,
3b~Naukova st., 79601 Lviv, Ukraine and Lviv National University,
1 Universytetska st.,
79602 Lviv, Ukraine}%
\email{rhryniv@iapmm.lviv.ua, yamykytyuk@yahoo.com}%

\address{\emph{Current address of R.H.:}
Institut f\"ur Angewandte Mathematik, Abteilung f\"ur
Wahr\-schein\-lich\-keitstheorie und Mathematische Statistik,
Wegelerstr.~6, D-53115 Bonn, Germany}
\email{rhryniv@wiener.iam.uni-bonn.de}

\thanks{${}^{\dag}$The work was partially supported by Ukrainian Foundation
for Basic Research, grant No.~01.07/00172.
R.~H. acknowledges the support of the Alexander von Humboldt foundation}%

\subjclass[2000]{Primary 34A55, Secondary 34B24, 34L05, 34L20}%
\keywords{Inverse spectral problems, Sturm-Liouville operators, singular
potentials, Sobolev spaces}%

\date{\today}%

%%%%%%%%%%%%%%%%%%%%%%%%%%%%%%%%%%%%%%%%%%%
\begin{abstract}
We solve the inverse spectral problems for the class of
Sturm--Liouville operators with singular real-valued potentials
from the Sobolev space $W^{s-1}_2(0,1)$, $s\in[0,1]$. The
potential is recovered from two spectra or from the spectrum and
norming constants. Necessary and sufficient conditions on the
spectral data to correspond to the potential in~$W^{s-1}_2(0,1)$
are established.
\end{abstract}
\maketitle

%%%%%%%%%%%%%%%%%%%%%%%%%%%%%%%%%%%%%%%%%%%
\section{Introduction}\label{sec:intr}

Suppose that $q$ is a real-valued distribution from~$W^{-1}_2(0,1)$. We
denote by $(\la_n^2)$ and $(\mu_n^2)$, $n\in\bN$, eigenvalues of
Sturm--Liouville operators $T_{\mathrm D}$ and $T_{\mathrm N}$
generated by the differential expression $-\frac{d^2}{dx^2} +q$ subject
to the Dirichlet and Neumann--Dirichlet boundary conditions
respectively (see Section~\ref{sec:pre} for precise definitions). The
operators $T_{\mathrm D}$ and $T_{\mathrm N}$ are bounded
below~\cite{SS} and thus become positive after addition of a suitable
constant to the potential~$q$. Henceforth there is no loss of
generality in assuming that the operators $T_{\mathrm D}$ and
$T_{\mathrm N}$ (and the numbers $\la_n$ and $\mu_n$) are positive. The
eigenvalues $\la_n^2$ and $\mu_n^2$ are simple; we arrange them in
increasing order and recall the following their
properties~\cite{An,HMtwo,SS1}:
\begin{itemize}
\item [(A1)] the sequences $(\la_n)$ and $(\mu_n)$ interlace, i.e.,
    $\mu_n<\la_n<\mu_{n+1}$ for all $n\in\bN$;
\item [(A2)] the numbers $\la_n$ and $\mu_n$ obey the asymptotics
    \begin{equation}\label{eq:lamu}
        \la_n = \pi n + \wt\la_n, \qquad \mu_n=\pi(n-1/2) + \wt\mu_n
    \end{equation}
 with some $\ell_2$-sequences $(\wt\la_n)$ and $(\wt\mu_n)$.
\end{itemize}

Conversely, it was shown in~\cite{HMtwo} that if two sequences
$(\la_n^2)$ and $(\mu_n^2)$ of positive numbers satisfy properties
(A1) and (A2), then there exists a unique $q\in W^{-1}_2(0,1)$
such that $(\la_n^2)$ and $(\mu_n^2)$ are eigenvalues of the
Sturm--Liouville operators $T_{\mathrm D}$ and $T_{\mathrm N}$
with potential~$q$. In other words, the inverse spectral problem
of recovering the potential by two spectra is uniquely soluble in
the class of Sturm--Liouville operators with singular potentials
from~$W^{-1}_2(0,1)$. The papers~\cite{HMinv,HMtwo} give the
corresponding reconstruction algorithm and thus extend the
classical inverse spectral theory for Sturm--Liouville operators
developed by Gelfand, Levitan, Marchenko, and Krein in
1950-ies~\cite{GL,K,M}.

The aim of the present work is to show that the above inverse spectral
problem is completely soluble in the class of Sturm--Liouville
operators with potentials from~$W^{s-1}_2(0,1)$ for every $s\in[0,1]$.
More precisely, we shall formulate necessary and sufficient conditions
on sequences $(\la_n^2)$ and $(\mu_n^2)$ in order that they are
Dirichlet and Neumann--Dirichlet eigenvalues of a Sturm--Liouville
operator with singular potential~$q$ from~$W^{s-1}_2(0,1)$.

The particular case $s=1$ corresponds to potentials in $L_2(0,1)$;
the classical theorem by Marchenko~\cite[Theorem 3.4.1]{Ma} states
that necessary and sufficient conditions on the eigenvalues
$(\la_n^2)$ and $(\mu_n^2)$ are (A1) and (A2) with the
specification that $\wt\la_n$ and $\wt\mu_n$ have the form
\[
    \wt\la_n = \frac{c}n + \frac{\wh \la_n}{n}, \qquad
    \wt\mu_n = \frac{c}n + \frac{\wh \mu_n}{n}
\]
with $c\in\bR$ and some $\ell_2$-sequences $(\wh\la_n)$ and
$(\wh\mu_n)$.

For an arbitrary intermediate value $s\in(0,1)$, the direct spectral
problem was studied in~\cite{HMas,KM,SS1}. For instance, it was proved
in~\cite{HMas} that $\wt\la_n$ and $\wt\mu_n$ are respectively even and odd
sine Fourier coefficients of some function from~$W^{s}_2(0,1)$ (cf. the
above-mentioned cases $s=0$ and $s=1$). More exactly, the main result
from~\cite{HMas} reads as follows.

\begin{theorem*}{Theorem~A. }
Assume that $q\in W^{s-1}_2(0,1)$ for some $s\in[0,1]$ and that $\wt\la_n$,
$\wt\mu_n$ are defined through~\eqref{eq:lamu}. Then the function $\si^*$
given by
\begin{equation}\label{eq:si*}
    \si^*(x):=  2 \sum_{n=1}^\infty \wt\mu_n\sin[(2n-1)\pi x] -
                2 \sum_{n=1}^\infty \wt\la_n\sin (2\pi nx)
\end{equation}
belongs to~$W^s_2(0,1)$. Moreover, $\si^*-\si \in W^{2s}_2(0,1)$, where
$\si$ is any of the distributional primitives of~$q$.
\end{theorem*}

In the present paper we show that the condition $\si^*\in
W^s_2(0,1)$ is a sufficient addendum to (A1) and (A2) guaranteeing
that the corresponding potential belongs to~$W^{s-1}_2(0,1)$. Our
main result is as follows.

\begin{theorem}\label{thm:main1}
In order that two sequences $(\la_n^2)$ and $(\mu_n^2)$ be
eigenvalues of (positive) Sturm--Liouville operators $T_{\mathrm
D}$ and $T_{\mathrm N}$ with potential from $W^{s-1}_2(0,1)$,
$s\in[0,1]$, it is necessary and sufficient that assumptions (A1),
(A2) hold and that the function~$\si^*$ of~\eqref{eq:si*} belongs
to $W^s_2(0,1)$.
\end{theorem}

As an intermediate step we solve the inverse spectral problem of
recovering the potential of a Sturm--Liouville expression by its
Dirichlet spectrum $(\la_n^2)$ and the so-called norming
constants~$(\al_n)$. We recall that
\[
    \al_n:= \Bigl(2\int_0^1|u_n(x)|^2\,dx\Bigr)^{-1},
\]
where $u_n$ is an eigenfunction of the operator $T_{\mathrm D}$
that corresponds to the eigenvalue $\la_n^2$ and satisfies the
initial condition $u_n^{[1]}(0)=\la_n$, with $u^{[1]}$ denoting
the quasi-derivative of a function~$u$, see Section~\ref{sec:pre}.
Alternatively, we can reduce the inverse spectral problem by two
spectra to recovering the potential by the spectrum $(\mu_n^2)$ of
the operator $T_{\mathrm N}$ and the norming constants $(\be_n)$;
the latter are defined as
\[
    \be_n:= \Bigl(2\int_0^1 |v_n(x)|^2\,dx\Bigr)^{-1},
\]
where $v_n$ is an eigenfunction of the operator $T_{\mathrm N}$
that corresponds to the eigenvalue~$\mu_n^2$ and satisfies the
initial condition $v_n(0)=1$.

In the case $q\in W^{-1}_2(0,1)$ (i.e., for $s=0$) the norming constants
$\al_n$ and $\be_n$ have the asymptotics $\al_n = 1 + \wt\al_n$,
$\be_n=1+\wt\be_n$ with $\ell_2$-sequences $(\wt\al_n)$ and $(\wt\be_n)$,
see~\cite{HMinv}. For $s\in(0,1]$ this asymptotics refines as follows. We
introduce the function $\ga$ via the formula
\begin{equation}\label{eq:ga}
    \ga(x):=    2 \sum_{n=1}^\infty \wt\be_n \cos[(2n-1)\pi x] -
                2 \sum_{n=1}^\infty \wt\al_n \cos (2\pi nx)
\end{equation}
and also put
\[
    \ga^*(x) := -2x \si^*(1-x).
\]

\begin{theorem}\label{thm:albe}
Assume that $\si^*\in W^s_2(0,1)$ for some $s\in[0,1]$. Then the function
$\ga$ also belongs to $W^s_2(0,1)$; moreover, $\ga-\ga^*\in W^{2s}_2(0,1)$.
\end{theorem}

In particular, we see that if a primitive $\si$ of $q$ belongs to
$W^s_2(0,1)$, then the sequences $(-\wt\al_n)$ and $(\wt\be_n)$ are even
and odd cosine Fourier coefficients respectively of the function $\ga\in
W^s_2(0,1)$ given by~\eqref{eq:ga}. In the reverse direction the claim is
that if the even and odd parts of the functions $\si^*$ and $\ga$ belong to
$W^s_2(0,1)$, then $\si$ is an $W^s_2$-function (and hence the potential
$q$ belongs to~$W^{s-1}_2(0,1)$). More exactly, the following two
statements hold true.

\begin{theorem}\label{thm:main2}
Sequences $(\la_n^2)$ and $(\al_n)$ of positive numbers are
eigenvalues and norming constants for some Sturm--Liouville
operator $T_{\mathrm D}$ with real-valued potential $q\in
W^{s-1}_2(0,1)$, $s\in[0,1]$, if and only if the following
conditions are satisfied:
\begin{itemize}
\item [(B1)] the numbers $\la_1<\la_2<\dots$ obey the asymptotics
    $\la_n=\pi n + \wt\la_n$, where $\wt\la_n$ are even sine Fourier
    coefficients of some function from $W^s_2(0,1)$;
\item [(B2)] the numbers $\wt\al_n:=\al_n-1$ are even cosine Fourier
    coefficients of some function from $W^s_2(0,1)$.
\end{itemize}
\end{theorem}

\begin{theorem}\label{thm:main3}
Sequences $(\mu_n^2)$ and $(\be_n)$ of positive numbers are
eigenvalues and norming constants for some Sturm--Liouville
operator $T_{\mathrm N}$ with real-valued potential $q\in
W^{s-1}_2(0,1)$, $s\in[0,1]$ if and only if the following
conditions are satisfied:
\begin{itemize}
\item [(C1)] the numbers $\mu_1<\mu_2<\dots$ obey the asymptotics
    $\mu_n=\pi (n-1/2) + \wt\mu_n$, where $\wt\mu_n$ are odd sine Fourier
    coefficients of some function from $W^s_2(0,1)$;
\item [(C2)] the numbers $\wt\be_n:=\be_n-1$ are odd cosine Fourier
    coefficients of some function from $W^s_2(0,1)$.
\end{itemize}
\end{theorem}

The organization of the paper is as follows. In
Section~\ref{sec:pre} we give the precise definitions of the
operators $T_{\mathrm D}$ and $T_{\mathrm N}$. In
Section~\ref{sec:norm} asymptotics of the norming constants is
established, based on which Theorem~\ref{thm:albe} is proved. The
algorithm of solution of the inverse spectral problems under
consideration and the proofs of Theorems~\ref{thm:main1},
\ref{thm:main2}, and \ref{thm:main3} are given in
Section~\ref{sec:inverse}. Finally, in Appendix~\ref{sec:sobolev}
some necessary facts about Sobolev spaces $W_2^s(0,1)$ and Fourier
series therein are gathered.

%%%%%%%%%%%%%%%%%%%%%%%%%%%%%%%%%%%%%%%%%%%
\section{Preliminaries}\label{sec:pre}

Suppose that $q\in W^{-1}_2(0,1)$ is real-valued. We fix an arbitrary
real-valued distributional primitive $\si\in L_2(0,1)$ of~$q$ (so that
$q=\si'$ in the sense of distributions) and consider the differential
expression
\[
    l_\si(u) := -(u'-\si u)' -\si u'
\]
on its ``maximal'' domain in~$L_2(0,1)$,
\[
    \fD(l_\si) = \{ u \in W^1_1(0,1) \mid u^{[1]} \in W^1_1(0,1),\
        l_\si(u) \in L_2(0,1)\}.
\]
Here and hereafter, $u^{[1]}$ stands for the \emph{quasi-derivative}
$u'-\si u$ of a function $u$. It is easily seen that $l_\si(u) = -u'' + q
u$ in the sense of distributions, so that $l_\si$ is a
\emph{regularization} of the differential expression $-\frac{d^2}{dx^2}+q$.

The operators $T_{\mathrm D}$ and $T_{\mathrm N}$ are defined as
the restrictions of $l_\si$ imposing the corresponding boundary
conditions
\begin{align*}
    \fD(T_{\mathrm D}) &= \{ u \in \fD(l_\si) \mid u(0)=u(1)=0\},\\
    \fD(T_{\mathrm N}) &= \{ u \in \fD(l_\si) \mid u^{[1]}(0)=u(1)=0\}.
\end{align*}
It is known~\cite{SS} that the operators $T_{\mathrm D}$ and
$T_{\mathrm N}$ are selfadjoint, bounded below, and have discrete
spectra. In some sense they are the most natural Dirichlet and
Neumann--Dirichlet Sturm--Liouville operators associated with the
differential expression~$-\frac{d^2}{dx^2} + q$ with singular~$q$; see
the discussion in~\cite{SS1}.

%%%%%%%%%%%%%%%%%%%%%%%%%%%%%%%%%%%%%%%%%%%

\section{Asymptotics of the norming constants}\label{sec:norm}

Suppose that $\si\in L_2(0,1)$ is real-valued. We denote by
$u_\pm(\cdot,\la)$ and $v_-(\cdot,\la)$ solutions of the equation
 \(
    l_\si (u) = \la^2 u
 \)
that satisfy the initial conditions
\[
    u_-(0,\la) = u_+(1,\la) = v_-^{[1]}(0,\la) = 0, \qquad
    u_-^{[1]}(0,\la) = u_+^{[1]}(1,\la) = v_-(0,\la) = 1.
\]
Observe that according to the definition of~$l_\si$ the equation
$l_\si(u) = \la^2 u$ is to be regarded as the first-order system
\[
    \frac{d}{dx}\binom{u}{u^{[1]}} =
        \begin{pmatrix} \si & 1\\ -\si^2-\la^2 & -\si
        \end{pmatrix} \binom{u}{u^{[1]}}.
\]
Since the entries of the above matrix are integrable, this system enjoys
the standard existence and uniqueness properties; in particular, the
solutions $u_\pm(\cdot,\la)$ and $v_-(\cdot,\la)$ are well defined for all
$\la\in\bC$.

Set $\Phi(\la):=\la u_-(1,\la)$ and $\Psi(\la):=u_+^{[1]}(0,\la)$;
then the numbers $\pm\la_n$ and $\pm\mu_n$ are zeros of $\Phi$ and
$\Psi$ respectively. $\Phi$ and $\Psi$ are respectively odd and
even entire functions of order~$1$ and hence can be represented by
their canonical Hadamard products~\cite{HMtwo}, namely,
\[
    \Phi(\la) = \la\prod_{n=1}^\infty \frac{\la_n^2-\la^2}{\pi^2 n^2}, \qquad
    \Psi(\la) = \prod_{n=1}^\infty \frac{\mu_n^2-\la^2}{\pi^2 (n-1/2)^2}.
\]
Moreover, it turns out that the norming constants
$\al_n=(\sqrt2\la_n\|u_-(\cdot,\la_n)\|)^{-2}$ and
$\be_n=(\sqrt2\|v_-(\cdot,\mu_n)\|)^{-2}$ can be expressed via the
functions~$\Phi$ and $\Psi$ only (cf.~\cite{GS,HMtwo}).

\begin{lemma}\label{lem:aln}
The norming constants~$\al_n$ and $\be_n$ satisfy the following equalities:
\begin{equation}\label{eq:alben}
    \al_n = \frac{\Psi(\la_n)}{\dot{\Phi}(\la_n)},\qquad
    \be_n =-\frac{\Phi(\mu_n)}{\dot{\Psi}(\mu_n)}.
\end{equation}
\end{lemma}

\begin{proof}
The Green function $G_{\mathrm D}(x,y,\la^2)$ of the
operator~$T_{\mathrm D}$ (i.e., the kernel of the resolvent
$(T_{\mathrm D}-\la^2)^{-1}$) equals
\[
    G_{\mathrm D}(x,y,\la^2) = \sum_{n=1}^{\infty}
        \frac{2\al_n\la_n^2 u_-(x,\la_n) u_-(y,\la_n)}{\la^2_n - \la^2}.
\]
On the other hand, we have
\[
    G_{\mathrm D}(x,y,\la^2) = \frac1{W(\la)} \left\{
        \begin{aligned}{u_-(x,\la) {u_+(y,\la)}},& \qquad
                        0\le x\le y\le1,\\
                    {u_-(y,\la) {u_+(x,\la)}},& \qquad
                        0\le y\le x\le1,
        \end{aligned} \right.
\]
where
     $W(\la):=u_+(x,\la)u^{[1]}_-(x,\la)
                -u_-(x,\la)u^{[1]}_+(x,\la)$
is the Wronskian of $u_+$ and $u_-$. The value of $W(\la)$ is
independent of $x\in[0,1]$; in particular, taking $x=1$ and $x=0$
we find that
\begin{equation}\label{eq:u+u-}
W(\la)\equiv -u_-(1,\la)\equiv u_+(0,\la).
\end{equation}

Equating the two expressions and comparing the residues at the poles
$\la=\la_n$, we find that
\[
    \al_n\la_n u_-(y,\la_n) =  \frac{u_+(y,\la_n)}{\dot{u}_-(1,\la_n)}.
\]
Observe that $u_-(\cdot,\la_n)$ and $u_+(\cdot,\la_n)$ are collinear and
hence
\[
    u_+(\cdot,\la_n)/u_-(\cdot,\la_n) =
    u_+^{[1]}(\cdot,\la_n)/u_-^{[1]}(\cdot,\la_n) =
    u_+^{[1]}(0,\la_n)=\Psi(\la_n);
\]
combining the above relations, we conclude that
\[
    \al_n = \frac{u_+(y,\la_n)}{u_-(y,\la_n)}
            \frac{1}{\la_n\dot{u}_-(1,\la_n)}
          = \frac{\Psi(\la_n)}{\dot{\Phi}(\la_n)}
\]
as claimed.

In a similar fashion we equate two expressions for the Green's
function~$G_{\mathrm N}(x,y,\la^2)$ of the operator $T_{\mathrm
N}$, namely,
\[
    \sum_{n=1}^\infty \frac{2\be_n v_-(x,\mu_n) v_-(y,\mu_n)}{\mu_n^2-\la^2}
    \equiv
    \frac1{W_1(\la)} \left\{
        \begin{aligned}{v_-(x,\la) {u_+(y,\la)}},& \qquad
                        0\le x\le y\le1,\\
                    {v_-(y,\la) {u_+(x,\la)}},& \qquad
                        0\le y\le x\le1.
        \end{aligned} \right.
\]
Here $W_1$ is the Wronskian of $u_+$ and $v_-$ and it is identically equal
to $-\Psi$, as follows from the equalities
\[
    W_1(\la):=u_+(x,\la)v^{[1]}_-(x,\la)-v_-(x,\la)u^{[1]}_+(x,\la)
        = -u^{[1]}_+(0,\la) = -\Psi(\la).
\]
Therefore we find that
\[
    \be_n =\frac{\mu_n u_+(y,\mu_n)}{v_-(y,\mu_n)\dot{\Psi}(\mu_n)}
        = \frac{\mu_n u_+(0,\mu_n)}{\dot{\Psi}(\mu_n)}
        =-\frac{\mu_n u_-(1,\mu_n)}{\dot{\Psi}(\mu_n)}
        =-\frac{\Phi(\mu_n)}{\dot{\Psi}(\mu_n)},
\]
where the second equality is obtained by taking $y=0$, while the third one
follows from~\eqref{eq:u+u-}. The lemma is proved.
\end{proof}

In what follows, we shall say that a function $f\in L_2(0,1)$ is \emph{odd}
(respectively \emph{even\/}) if $f(1-x)\equiv -f(x)$ (respectively, if
$f(1-x)\equiv f(x)$). Denote by $L_{2,\mathrm{o}}(0,1)$ and
$L_{2,\mathrm{e}}(0,1)$ the subspaces of $L_2(0,1)$ consisting of odd and
even functions respectively. We shall denote by $f_{\mathrm{o}}$ and
$f_{\mathrm{e}}$ respectively the odd and even parts of a function $f$;
obviously,
\[
    f_{\mathrm{o}}(x) = \tfrac12 \bigl[f(x)-f(1-x)\bigr], \qquad
    f_{\mathrm{e}}(x) = \tfrac12 \bigl[f(x)+f(1-x)\bigr].
\]

\begin{lemma}\label{lem:int}
The functions $\Phi$ and $\Psi$ admit the integral representations
\begin{equation}\label{eq:int}
\begin{aligned}
    \Phi(\la) &= \sin\la + \int_0^1 \phi(x) \sin [\la (1-2x)]\,dx,\\
    \Psi(\la) &= \cos\la + \int_0^1 \psi(x) \cos [\la (1-2x)]\,dx,
\end{aligned}
\end{equation}
in which $\phi\in L_{2,\mathrm{o}}(0,1)$ and $\psi\in
L_{2,\mathrm{e}}(0,1)$.
\end{lemma}

\begin{proof}
Using the technique of the transformation operators~\cite{HMtr}, $\Phi$ and
$\Psi$ can be shown to admit the integral representations of the form
\begin{align*}
    \Phi(\la) &= \sin\la + \int_0^1 \wt\phi(x) \sin \la x\,dx,\\
    \Psi(\la) &= \cos\la + \int_0^1 \wt\psi(x) \cos \la x\,dx
\end{align*}
with some $L_2$-functions $\wt\phi$ and $\wt\psi$; see detailed
derivation in~\cite{HMtwo}. Now we put
\begin{align*}
    \phi(x) &:= \begin{cases}
                 \wt\phi(1-2x) \quad &\text{if}\quad x\in[0,1/2]\\
                -\wt\phi(2x-1) \quad &\text{if}\quad x\in(1/2,1]
                \end{cases},\\
    \psi(x) &:= \begin{cases}
                 \wt\psi(1-2x) \qquad &\text{if}\quad x\in[0,1/2]\\
                 \wt\psi(2x-1) \qquad &\text{if}\quad x\in(1/2,1].
                \end{cases}
\end{align*}
It is easily seen that $\phi\in L_{2,\mathrm{o}}(0,1)$, $\psi\in
L_{2,\mathrm{e}}(0,1)$, and that equalities~\eqref{eq:int} hold. The lemma
is proved.
\end{proof}

The next lemma tells us that the values of $\Phi$ and $\Psi$ at the points
$\la_n$ and $\mu_n$ are expressed through sine and cosine Fourier
coefficients of some related functions. In the following $s_n(f)$ and
$c_n(f)$ will stand for respectively $n$-th sine and $n$-th cosine Fourier
coefficients of a function $f\in L_2(0,1)$; see~\eqref{eq:frtr} for exact
formulae. We also denote by $S$ the operator of multiplication by $1-2x$,
i.e., $(Sf)(x) = (1-2x)f(x)$.

\begin{lemma}\label{lem:sncn}
For an arbitrary $f\in L_2(0,1)$, the following equalities hold:
\begin{enumerate}
\item $\int_0^1 f(x) \sin[\la_n(1-2x)]\,dx =
            (-1)^{n+1}\bigl[s_{2n}(f) - \wt\la_n c_{2n}(Sf) +
            \wt\la_n^2 s_{2n}(f_1)\bigr]$;
\item $\int_0^1 f(x) \cos[\la_n(1-2x)]\,dx =
            (-1)^n\bigl[c_{2n}(f) + \wt\la_n s_{2n}(Sf) +
            \wt\la_n^2 c_{2n}(f_2)\bigr]$;
\item $\int_0^1 f(x) \sin[\mu_n(1-2x)]\,dx =
            (-1)^{n+1}\bigl[c_{2n-1}(f) + \wt\mu_n s_{2n-1}(Sf) +
            \wt\mu_n^2 c_{2n-1}(f_3)\bigr]$;
\item $\int_0^1 f(x) \cos[\mu_n(1-2x)]\,dx =
            (-1)^{n+1}\bigl[s_{2n-1}(f) - \wt\mu_n c_{2n-1}(Sf) +
            \wt\mu_n^2 s_{2n-1}(f_4)\bigr]$,
\end{enumerate}
where $f_j$, $j=1,2,3,4$, are some functions from $L_2(0,1)$.
\end{lemma}

\begin{proof}
We shall prove only part (1) as the other parts are established
analogously. Using the equality
\begin{align*}
    \sin [\la_n (1-2x)] &= (-1)^{n+1} \cos[\wt\la_n(1-2x)]\sin (2\pi nx)\\
        &\quad + (-1)^n  \sin[\wt\la_n(1-2x)]\cos (2\pi nx),
\end{align*}
the asymptotic relations
\[
    \sin t = t + O(t^3), \qquad \cos t = 1- t^2/2 + O(t^4), \quad t\to0,
\]
and the fact that $(\wt\la_n)\in\ell_2$, we find that
\[
    \int_0^1 f(x)\sin [\la_n (1-2x)]\,dx =
        (-1)^{n+1} \bigl[ s_{2n}(f) - \wt\la_n c_{2n}(Sf)
            + \wt\la_n^2 a_n\bigr]
\]
for some $\ell_2$-sequence $(a_n)$. Clearly, there exists a function
$f_1\in L_2(0,1)$ such that $a_n = s_{2n}(f_1)$ for all $n\in\bN$ and the
proof of part (1) is complete.
\end{proof}

\begin{remark}\label{rem:sncn}
Put
\begin{alignat*}{2}
    g_1&:= \si^* \tast Sf  + \si^* \tast(\si^* \hast f_1),&\qquad
    g_2&:= -\si^* \hast Sf + \si^* \hast(\si^* \tast f_2),\\
    g_3&:= \si^* \hast Sf  + \si^* \hast(\si^* \tast f_3),&\qquad
    g_4&:= -\si^* \tast Sf + \si^* \tast(\si^* \hast f_4) ,
\end{alignat*}
where the operations $\hast$ and $\tast$ are introduced in
Appendix~\ref{sec:sobolev}. By virtue of Lemma~\ref{lem:conv} we can
restate equalities (1)--(4) of the previous lemma as follows:
\begin{itemize}
\item[(1$'$)] $\int_0^1 f(x) \sin[\la_n(1-2x)]\,dx =
            (-1)^{n+1} s_{2n}(f+g_1) $;
\item[(2$'$)] $\int_0^1 f(x) \cos[\la_n(1-2x)]\,dx =
            (-1)^n c_{2n}(f+g_2)$;
\item[(3$'$)] $\int_0^1 f(x) \sin[\mu_n(1-2x)]\,dx =
            (-1)^{n+1} c_{2n-1}(f+g_3)$;
\item[(4$'$)] $\int_0^1 f(x) \cos[\mu_n(1-2x)]\,dx =
            (-1)^{n+1} s_{2n-1}(f+g_4)$.
\end{itemize}
If the functions $\si^*$ and $f$ belong to $W^s_2(0,1)$ for some
$s\in[0,1]$, then $Sf\in W^s_2(0,1)$ by Proposition~\ref{pro:interp}, and
thus Corollary~\ref{cor:conv} implies that the above functions $g_j$,
$j=1,2,3,4$, belong to $W^{2s}_2(0,1)$.
\end{remark}

Using~\eqref{eq:alben}, integral representations for $\Phi$ and $\Psi$, and
asymptotics of $\la_n$ and $\mu_n$, we can show that the norming constants
$\al_n$ and $\be_n$ obey the asymptotics $\al_n = 1 + \wt\al_n$ and $\be_n
= 1 +\wt\be_n$ with $\ell_2$-sequences $(\wt\al_n)$ and $(\wt\be_n)$. It
turns out that if the spectral data $(\la_n^2)$ and $(\mu_n^2)$ have better
asymptotics, then the functions $\phi$ and $\psi$ in~\eqref{eq:int} become
smoother, and the asymptotics of $\al_n$ and $\be_n$ refine.

\begin{lemma}\label{lem:phipsi} Assume that the numbers
$\wt\la_n:=\la_n-\pi n$ and $\wt\mu_n =\mu_n-\pi (n-1/2)$ are such that the
function $\si^*$ of~\eqref{eq:si*} belongs to $W^s_2(0,1)$ for some
$s\in[0,1]$. Then the functions $\phi$ and $\psi$ in integral
representation~\eqref{eq:int} of~$\Phi$ and~$\Psi$ have the form
\[
    \phi =-\si^*_{\mathrm{o}} + \phi_1, \qquad
    \psi = \si^*_{\mathrm{e}} + \psi_1,
\]
where $\phi_1$ and $\psi_1$ are respectively some odd and even functions
from $W^{2s}_2(0,1)$.
\end{lemma}

\begin{proof}
In virtue of Lemma~\ref{lem:int} the equality $\Phi(\la_n)=0$ can be recast
as
\[
    \sin\la_n + \int_0^1 \phi(x) \sin[\la_n(1-2x)]\,dx = 0.
\]
Observe that $\sin\la_n=(-1)^n\sin\wt\la_n$ and that $\sin\wt\la_n =
\wt\la_n + \wt\la_n^2 b_n$ for some $\ell_2$-sequence~$(b_n)$. Combining
this observation with Lemma~\ref{lem:sncn}, we arrive at the relation
\[
     \wt\la_n - s_{2n}(\phi) + \wt\la_n c_{2n}(S\phi)
                       +\wt\la_n^2 s_{2n}(\wh\phi)=0
\]
for some odd function $\wh\phi\in L_{2,\mathrm{o}}(0,1)$. Using
Lemma~\ref{lem:conv} and recalling that $\wt\la_n=-s_{2n}(\si^*) =
-s_{2n}(\si^*_{\mathrm{o}})$, we conclude that
\[
    \phi = -\si^*_{\mathrm{o}} - \si^*_{\mathrm{o}} \tast
    \bigl[(S\phi) - \si^*_{\mathrm{o}} \hast \wh\phi \bigr].
\]
In particular, $\phi\in W^s_2(0,1)$ by Corollary~\ref{cor:conv}, so that
the function $(S\phi) - \si^*_{\mathrm{o}} \hast \wh\phi$ belongs to
$W^s_2(0,1)$, and again by Corollary~\ref{cor:conv} we get
$\phi_1:=\phi+\si^*_{\mathrm{o}} = \si^*_{\mathrm{o}} \tast
    \bigl[\si^*_{\mathrm{o}} \hast \wh\phi -(S\phi)\bigr]\in W^{2s}_2(0,1)$.
The fact that $\phi_1\in L_{2,\mathrm{o}}(0,1)$ is obvious.

In a similar manner, using the relations $\Psi(\mu_n)=0$ and
$\mu_n = \pi (n-\tfrac12) + \wt\mu_n$ and Lemmata~\ref{lem:int}
and \ref{lem:sncn} we find that
\[
    \wt\mu_n - s_{2n-1}(\psi) + \wt\mu_n c_{2n-1}(S\psi) +
                \wt\mu_n^2 s_{2n-1}(\wh\psi) =0
\]
for some function $\wh\psi\in L_{2,\mathrm{e}}(0,1)$. Replicating the
above reasoning, we conclude that $\psi = \si^*_{\mathrm{e}} + \psi_1$
for some even function $\psi_1$ from $W^{2s}_2(0,1)$ as claimed. The
proof is complete.
\end{proof}

\begin{proof}[Proof of Theorem~\ref{thm:albe}]
Formula~\eqref{eq:alben} implies that
\[
    \wt\al_n=\al_n-1 =
        \frac{\Psi(\la_n)-\dot{\Phi}(\la_n)}{\dot{\Phi}(\la_n)}, \qquad
    \wt\be_n=\be_n-1 =
        \frac{-\Phi(\mu_n)-\dot{\Psi}(\mu_n)}{\dot{\Psi}(\mu_n)}
\]
for all $n\in\bN$. According to Lemma~\ref{lem:int} we have
\begin{align*}
    \Psi(\la_n)-\dot{\Phi}(\la_n) &
            = \int_0^1 \th_1(x)\cos[\la_n(1-2x)]\,dx,\\
    -\Phi(\mu_n)-\dot{\Psi}(\mu_n) &
            = \int_0^1 \th_2(x)\sin[\mu_n(1-2x)]\,dx
\end{align*}
with $\th_1:= \psi - S\phi$ and $\th_2:=-\phi + S\psi$. Further, in view of
Lemma~\ref{lem:sncn} and Remark~\ref{rem:sncn},
\begin{align*}
    \int_0^1 \th_1(x) \cos[\la_n(1-2x)]\,dx
            &= (-1)^n c_{2n}(\th_1 + \wt\th_1),\\
    \int_0^1 \th_2(x) \sin[\mu_n(1-2x)]\,dx
            &= (-1)^{n+1} c_{2n-1}(\th_2 + \wt\th_2)
\end{align*}
with some functions $\wt\th_1$ and $\wt\th_2$ from $W^{2s}_2(0,1)$.

It follows from Lemma~\ref{lem:int} that
\begin{align*}
    \dot{\Phi}(\la_n) &= \cos\la_n
        + \int_0^1 (1-2x)\phi(x) \cos [\la_n (1-2x)]\,dx,\\
   -\dot{\Psi}(\mu_n) &= \sin\mu_n
        + \int_0^1 (1-2x)\psi(x) \sin [\mu_n (1-2x)]\,dx.
\end{align*}
We have $\cos\la_n = (-1)^n\cos\wt\la_n = (-1)^n (1 + \wt\la_n d_n)$ and
$\sin\mu_n = (-1)^{n+1} \cos\wt\mu_n =(-1)^{n+1} (1+ \wt\mu_n e_n)$ for
some $\ell_2$-sequences $(d_n)$ and $(e_n)$. Using Lemma~\ref{lem:sncn} and
Remark~\ref{rem:sncn}, we now conclude that
\[
    (-1)^n \dot{\Phi}(\la_n) = 1 + c_{2n}(g_1), \qquad
    (-1)^{n+1} \dot{\Psi}(\mu_n) = 1 + c_{2n-1}(g_2)
\]
for some functions $g_1$ and $g_2$ from $W^s_2(0,1)$. Since
$\dot{\Phi}(\la_n)\ne0$ and $\dot{\Psi}(\mu_n)\ne0$ for all $n\in\bN$,
Lemma~\ref{lem:wiener} implies that
\[
    \frac{(-1)^n}{\dot{\Phi}(\la_n)}     = 1 + c_{2n}(h_1), \qquad
    \frac{(-1)^{n+1}}{\dot{\Psi}(\mu_n)} = 1 + c_{2n-1}(h_2)
\]
for some functions $h_1$ and $h_2$ from $W^s_2(0,1)$.

We now combine the above relations to conclude that
\[
    \wt\al_n = c_{2n}(\th_1 + \wt\th_1)\bigl(1 + c_{2n}(h_1)\bigr), \qquad
    \wt\be_n = c_{2n-1}(\th_2 + \wt\th_2)\bigl(1 + c_{2n-1}(h_2)\bigr).
\]
It follows that $\ga = -\th_1+\th_2 + \wt\th$ for some $\wt\th\in
W^{2s}_2(0,1)$. Since by Lemma~\ref{lem:phipsi}
\begin{align*}
    -\th_1 + \th_2 = (S-I)(\phi+\psi)
        &= (S-I)(\si^*_{\mathrm{e}} -
            \si^*_{\mathrm{o}}) + (S-I)(\phi_1 + \psi_1)\\
        &= \ga^* + (S-I)(\phi_1 + \psi_1)
\end{align*}
and $(S-I)(\phi_1 + \psi_1)\in W^{2s}_2(0,1)$ by
Proposition~\ref{pro:interp}, we conclude that the function $\ga-\ga^*$ is
in $W^{2s}_2(0,1)$ as required. The theorem is proved.
\end{proof}

Observe that Theorem~\ref{thm:albe} gives necessary parts of
Theorems~\ref{thm:main2} and \ref{thm:main3}. Sufficient parts of these
theorems constitute the inverse spectral problem and are treated in the
next section.

%%%%%%%%%%%%%%%%%%%%%%%%%%%%%%%%%%%%%%%%%%%

\section{The inverse problem}\label{sec:inverse}

We start by recalling briefly the standard method of recovering
the potential of a Sturm--Liouville operator from the spectral
data---sequences of eigenvalues and the corresponding norming
constants. This method was suggested by Gelfand and Levitan
in~\cite{GL} for the case of regular (i.e., locally integrable)
potentials and was further developed in~\cite{HMinv} to cover
singular potentials from $W^{-1}_2(0,1)$.

Consider the functions
\begin{alignat*}{2}
    \om_1(x) &:= \sum_{n=1}^\infty
        \bigl[ \al_n \cos \la_n x - \cos(\pi nx)\bigr],& \qquad
                &x\in[0,2];\\
    \om_2(x) &:= \sum_{n=1}^\infty
        \bigl\{\be_n \cos \mu_n x - \cos[\pi (n-\tfrac12)x]\bigr\},& \qquad
                &x\in[0,2],
\end{alignat*}
which belong to $L_2(0,2)$ as soon as the functions $\si^*$ and $\ga$ are
in $L_2(0,1)$ (cf.~Lemmata~\ref{lem:om1} and \ref{lem:om2} below) and put
for $x,y\in[0,1]$
\begin{equation}\label{eq:f}
    f_j(x,y) := \om_j(|x-y|) + (-1)^j \om_j(x+y), \qquad j=1,2.
\end{equation}
We introduce an integral operator $F_j$ with kernel $f_j$; namely, $F_j$
acts in $L_2(0,1)$ according to the formula
\begin{equation}\label{eq:F}
    (F_j u)(x) := \int_0^1 f_j(x,y) u(y)\,dy.
\end{equation}

Let also $K_1$ and $K_2$ be the transformation operators for
$T_{\mathrm D}$ and $T_{\mathrm N}$ respectively. Recall that
$K_j$, $j=1,2$, is an integral operator with lower-triangular
kernel $k_j$, i.e., $k_j(x,y)=0$ a.e.~on the set $\{(x,y) \mid 0
\le x < t \le1\}$ and thus
\[
    (K_j u)(x) = \int_0^x k_j(x,y) u(y)\,dy, \qquad j=1,2.
\]
The operator $I+K_1$ transforms solutions of the unperturbed
equation $l_{0}(u)=\la^2 u$ (i.e., corresponding to $\si\equiv0$)
subject to the Dirichlet initial condition $u(0)=0$ into the
solutions of the equation $l_\si(u)=\la^2u$ subject to the
Dirichlet initial condition; the operator $I+K_2$ does the same
for the Neumann boundary condition at $x=0$.

Moreover, $f_j$ and $k_j$ are related through the so-called
Gelfand--Levitan--Marchenko (GLM) equation
\begin{equation}\label{eq:GLM}
    f_j(x,y) + k_j(x,y) + \int_0^x k_j(x,t)f_j(t,y)\,dt =0, \qquad
        0 \le y\le x \le 1.
\end{equation}
It is known~\cite{GKvolt,HMinv} that this GLM equation is naturally
related to the problem of factorisation of the operator $I+F_j$ in a
special manner and that~\eqref{eq:GLM} is uniquely soluble for~$k_j$ as
soon as the operator $I+F_j$ is (uniformly) positive. We note that
under properties (B1) and (B2) with $s=0$ the operator $I+F_1$ is
uniformly positive in $L_2(0,1)$, and the same is true of $I+F_2$ if
(C1) and (C2) hold with $s=0$ (see the details in~\cite{HMinv}).
Henceforth in both cases of interest the GLM equation possesses a
unique solution. Moreover, up to an additive constant $C_j$, the
primitive $\si$ of the potential $q$ equals
\begin{equation}\label{eq:si}
    \si(x) = (-1)^{j-1}2 \om_j(2x) -
        2\int_0^x k_j(x,t) f_j(t,x)\,dt + C_j.
\end{equation}
In what follows, we shall restrict ourselves to the case of the
operator $T_{\mathrm D}$; the proofs for the operator $T_{\mathrm
N}$ require only minor modifications.

In order to prove sufficiency part of Theorem~\ref{thm:main2} we need to
show first that the function $\om_1$ belongs to $W^s_2(0,1)$, then
establish some properties of the kernel $k_1$, and finally use
formula~\eqref{eq:si} to prove the inclusion $\si\in W^s_2(0,1)$.

\begin{lemma}\label{lem:om1}
Assume that the numbers $\wt\la_n$ and $\wt\al_n>-1$ are such that
there exist functions $g$ and $h$ in $W^s_2(0,1)$ with the property
that $\wt\la_n = s_{2n}(g)$ and $\wt\al_n=c_{2n}(h)$. Then the
function~$\om_1$ belongs to $W^s_2(0,2)$.
\end{lemma}

\begin{proof}
Observe first that, by the construction of $\si^*$ and $\ga$, we have
$\si^*_{\mathrm{o}}= g_{\mathrm{o}}$ and
$\ga_{\mathrm{e}}=h_{\mathrm{e}}$, whence $\si^*_{\mathrm{o}}\in
W^s_2(0,1)$ and $\ga_{\mathrm{e}}\in W^s_2(0,1)$. We write
\begin{align*}
    2\om_1(2x) &= 2\sum_{n=1}^\infty
        \bigl[(1+\wt\al_n)\cos(2\pi nx+2\wt\la_nx) - \cos(2\pi nx)\bigr]\\
    &= 2\sum_{n=1}^\infty \wt\al_n \cos(2\pi nx)
       + 2\sum_{n=1}^\infty (1+\wt\al_n)
                [\cos(2\wt\la_nx)-1] \cos(2\pi nx)\\
       & \quad -2\sum_{n=1}^\infty (1+\wt\al_n)
                \sin(2\wt\la_nx) \sin(2\pi nx)\\
        &=: -\ga_{\mathrm{e}}(x) + g_1(x) -  g_2(x),
\end{align*}
so that it remains to prove that the functions $g_1$ and $g_2$ belong
to $W^{s}_2(0,1)$.

Justification of the inclusions $g_1\in W^{s}_2(0,1)$ and $g_2\in
W^{s}_2(0,1)$ is similar, and we shall give it in detail only for the
function $g_1$. We have
\[
    g_1(x)=2\sum_{n=1}^\infty (1+\wt\al_n) \cos(2\pi nx)
        \sum_{k=1}^\infty (-1)^k\frac{(2\wt\la_nx)^{2k}}{(2k)!}.
\]
For $x\in[0,1]$, the estimate
\[
    \sum_{k=1}^\infty \frac{(2\wt\la_nx)^{2k}}{(2k)!}
        \le \cosh(2\wt\la_n)- 1 = \mathrm{O}(\wt\la_n^2)
\]
and the inclusion $(\wt\la_n)\in\ell_2$ imply that the above double
series for $g_1$ converges uniformly and absolutely. Changing the
summation order, we find that
\begin{equation}\label{eq:g1}
    g_1(x) = \sum_{k=1}^\infty \frac{(-1)^k2^{2k}}{(2k)!}x^{2k} h_k(x),
\end{equation}
where
\[
    h_k(x) := 2\sum_{n=1}^\infty (1+\wt\al_n)\,\wt\la_n^{2k}\cos(2\pi nx).
\]

Put $\tau:= \si^*_{\mathrm{o}}\tast \si^*_{\mathrm{o}}$; then by virtue
of Lemma~\ref{lem:conv} we find that
\[
    h_k = (\underbrace{\tau \oast \tau \oast \cdots
                    \oast \tau}_{k\ \text{times}}) +
        \ga_{\mathrm{e}} \oast (\underbrace{\tau \oast \tau \oast \cdots
                    \oast \tau}_{k\ \text{times}});
\]
see the definition of $\oast$ and $\tast$ in
Appendix~\ref{sec:sobolev}. Corollary~\ref{cor:conv} yields the
inclusion $h_k\in W^{s}_2(0,1)$; moreover, there exists a number
$C_1>0$ such that
\[
    \|h_k\|_{s} \le C_1 \bigl(1+\|\ga_\mathrm{e}\|_{0}\bigr)
        \ \|\si^*_{\mathrm{o}}\|^{2k}_{s}.
\]

Denote by $V$ the operator of multiplication by $x$; then $V^{2k}$ is
bounded both in $L_2(0,1)$ and in $W^1_2(0,1)$ and there exists
$C_2\ge1$ such that, for all $f\in W^1_2(0,1)$,
\[
    \|V^{2k}f\|_{0}\le \|f\|_{0}, \qquad
    \|V^{2k}f\|_{1} \le C_2 k \|f\|_{1}
\]
Interpolating between $W^1_2(0,1)$ and $L_2(0,1)$, we conclude that
$V^{2k}$ is bounded in every intermediate space $W^{s}_2(0,1)$ and
$\|V^{2k}f\|_{s} \le C_2 k^{s}\|f\|_{s}$ for all $f\in W^{s}_2(0,1)$.
Combining the above relations, we conclude that the series
in~\eqref{eq:g1} converges in~$W^{s}_2(0,1)$. Henceforth $g_1\in
W^{s}_2(0,1)$, which completes the proof.
\end{proof}

Denote by $\fA_s$ the set of all integral operators $K$ over~$(0,1)$,
whose kernels $k$ possess the following properties:
\begin{itemize}
\item [(1)] for every $x\in [0,1]$ the functions $k(x,\cdot)$ and
    $k(\cdot,x)$ belong to $W^s_2(0,1)$;
\item [(2)] the mappings
\[
    [0,1] \ni x \mapsto k(x,\cdot) \in W^s_2(0,1), \qquad
    [0,1] \ni x \mapsto k(\cdot,x) \in W^s_2(0,1)
\]
are continuous.
\end{itemize}

The results of~\cite{My} imply the following statement.

\begin{proposition}\label{pro:fact}
Assume that $F$ is an integral operator with kernel $f$ such that
$F\in\fA_s$, $s\in[0,\frac12)$, and $I+F>0$. Let $k$ be equal to a
solution of the corresponding GLM equation~\eqref{eq:GLM} in the domain
$0\le y\le x\le1$ and be zero in the domain $0\le x<y\le1$; then an
integral operator~$K$ with kernel $k$ also belongs to~$\fA_s$.
\end{proposition}

Let the assumptions of Lemma~\ref{lem:om1} hold for some $s\in
[0,\frac12)$, so that $\om_1\in W^s_2(0,1)$. Then the operator $F_1$ given
by~\eqref{eq:f}--\eqref{eq:F} belongs to $\fA_s$. Indeed, properties (1)
and (2) of the definition of $\fA_s$ for the kernel $f_1$ follow from the
fact that
\begin{itemize}
\item[(a)] the operator $P$ restricting a function on $\bR$
 onto $(0,1)$ is a bounded mapping from $W^s_2(\bR)$ into
 $W^s_2(0,1)$~\cite[Theorem~1.9.1]{LM};
\item[(b)] the translations $T_t f (\cdot) := f(\cdot+t)$,
    $t\in\bR$, form a $C_0$-group in $W^s_2(\bR)$ \cite{My}.
\end{itemize}

With these preliminaries in hand, we can complete the inverse spectral
analysis of Theorems~\ref{thm:main1}, \ref{thm:main2}, and \ref{thm:main3}.

\begin{proof}[Proof of Theorem~\ref{thm:main1}.]
The necessity part of the theorem follows from Theorem~A, hence we need to
prove only the sufficiency part, i.e., that properties (A1) and (A2) and
the inclusion $\si^*\in W^s_2(0,1)$ imply that $\si\in W^s_2(0,1)$.

Assume therefore that the sequences $(\la_n^2)$ and $(\mu_n^2)$ satisfy
properties (A1) and (A2). Applying the reconstruction procedure
explained above (and developed in detail for the case of singular
potentials from~$W^{-1}_2(0,1)$ in~\cite{HMtwo}), we find a unique
function $\si\in L_2(0,1)$ such that $\la_n^2$ and $\mu_n^2$ are
eigenvalues of the corresponding operators $T_{\mathrm D}$ and
$T_{\mathrm N}$. It remains to prove that the inclusion $\si^*\in
W^s_2(0,1)$ yields $\si\in W^s_2(0,1)$. We shall consider separately
the cases $s\in[0,\frac12)$, $s\in[\frac12,1)$, and $s=1$.

\textbf{Case 1: $s\in[0,\frac12)$.} By Theorem~\ref{thm:albe} the
function~$\ga$ belongs to~$W^s_2(0,1)$; hence by Lemma~\ref{lem:om1}
the function $\om_1$ is in $W^s_2(0,1)$ and, as explained above, the
integral operator~$F_1$ falls into the set~$\fA_s$. It follows from
Proposition~\ref{pro:fact} that the solution~$k_1$ of the GLM
equation~\eqref{eq:GLM} generates an integral operator $K_1$ that also
belongs to~$\fA_s$.

In view of formula~\eqref{eq:si} and the inclusion~$\om_1\in
W^s_2(0,1)$, the theorem will be proved as soon as we show that the
integral
 \(
    \int_0^x k_1(x,t) f_1(t,x)\,dt
 \)
defines a function from $W^s_2(0,1)$. Observe that $k_1(x,t)=0$
a.e.~for $0\le x < t\le1$, so that we can extend the range of
integration to $[0,1]$. Hence we put
\[
    \eta(x) := \int_0^1 k_1(x,t) f_1(t,x)\,dt
\]
and shall prove that $\eta \in W^s_2(0,1)$.

Recall \cite[Theorem~1.10.2]{LM} that one of the equivalent norms
in the space $W^s_2(0,1)$ is given by
\[
    \|\eta\|_{s} = \biggl( \|\eta\|^2_0
        + 2\int_0^1 \int_x^1
           \frac{|\eta(x)-\eta(y)|^2}{(x-y)^{1+2s}} \,dy\,dx\biggr)^{1/2}.
\]
Now we use the fact that the operators $K_1$ and $F_1$ belong to $\fA_s$.
In particular, there exists a constant $C_1>0$ such that
\[
    \max_{x\in [0,1]} \bigl( \|k_1(x,\cdot)\|^2_{s}+
        \|k_1(\cdot,x)\|^2_{s}+ \|f_1(x,\cdot)\|^2_{s}+
        \|f_1(\cdot,x)\|^2_{s}\bigr) \le C_1,
\]
and we use this inequality to derive the estimates
\[
    \|\eta\|^2_{0} \le \int_0^1 \Bigl(
        \int_0^1 |k_1(x,y)|^2\,dy \int_0^1 |f_1(y,x)|^2\,dy
            \Bigr) \,dx \le C_1^2
\]
and
\begin{align*}
    |\eta(x)-\eta(y)|^2 &\le
             2 \Bigl| \int_0^1 |k_1(x,t)-k_1(y,t)||f_1(t,x)|\,dt\Bigr|^2\\
    &\quad + 2 \Bigl| \int_0^1 |k_1(y,t)||f_1(t,x)-f_1(t,y)|\,dt\Bigr|^2\\
    &\le 2C_1 \Bigl( \int_0^1 |k_1(x,t)-k_1(y,t)|^2\,dt
        + \int_0^1 |f_1(t,x)-f_1(t,y)|^2\,dt\Bigr).
\end{align*}
It follows that
\[
    \|\eta\|^2_{s} \le C_1^2 +
        4C_1 \int_0^1 \Bigl[\|k_1(\cdot,t)\|^2_{s} +
            \|f_1(t,\cdot)\|^2_{s}\Bigr]\,dt \le 5C_1^2,
\]
so that $\eta\in W^s_2(0,1)$, and the proof for the case~$s\in[0,\frac12)$
is complete.

\textbf{Case 2: $s\in[\frac12,1)$.} Case 1 applied to the exponent
$\tfrac{s}2$ gives $\si\in W^{s/2}_2(0,1)$, so that $\si-\si^*\in
W^s_2(0,1)$ by Theorem~A. Since by assumption $\si^*\in W^s_2(0,1)$, we
have $\si\in W^s_2(0,1)$ as required.

\textbf{Case 3: $s=1$.} We again use the bootstrap method: first, by
Case~2, $\si\in W^t_2(0,1)$ for any $t\in[0,1)$, e.g., $t=\frac12$; then,
by Theorem~A, $\si-\si^*\in W^1_2(0,1)$. This inclusion yields $\si\in
W^1_2(0,1)$, and the proof is complete.
\end{proof}

\begin{proof}[Proof of Theorem~\ref{thm:main2}.]
If $\si\in W^s_2(0,1)$, then $\si^*\in W^s_2(0,1)$ by Theorem~A, hence
$\ga\in W^s_2(0,1)$ by Theorem~\ref{thm:albe}. Properties (B1) and (B2)
then obviously hold as, by the construction of $\si^*$ and $\ga$, $\wt\la_n
= s_{2n}(-\si^*)$ and $\wt\al_n = c_{2n}(-\ga)$.

Conversely, assume (B1) and (B2). According to the results
of~\cite{HMinv}, there exists a unique Sturm--Liouville operator
$T_{\mathrm D}$ with potential $q$ from $W^{-1}_2(0,1)$ such that
$\la_n^2$ and $\al_n$ are eigenvalues and norming constants of
$T_{\mathrm D}$. We need to prove that validity of (B1) and (B2)
implies that the recovered potential $q$ belongs in fact to
$W^{s-1}_2(0,1)$ (i.e., that any primitive $\si$ of $q$ belongs to
$W^s_2(0,1)$). In fact, under (B1) and (B2) Lemma~\ref{lem:om1} yields
$\om_1\in W^s_2(0,2)$, and it remains to observe that the proof of
Theorem~\ref{thm:main1} derives from this the inclusion $\si\in
W^s_2(0,1)$. The proof is complete.
\end{proof}

Proof of Theorem~\ref{thm:main3} is completely analogous; the only
essential reservation is that Lemma~\ref{lem:om1} should be replaced with
the following its counterpart (we leave both proofs to the reader):

\begin{lemma}\label{lem:om2}
Assume that the numbers $\wt\mu_n$ and $\wt\be_n>-1$ are such that
there exist functions $g$ and $h$ in $W^s_2(0,1)$ with the property
that $\wt\mu_n = s_{2n-1}(g)$ and $\wt\be_n=c_{2n-1}(h)$. Then the
function~$\om_2$ belongs to $W^s_2(0,1)$.
\end{lemma}

%%%%%%%%%%%%%%%%%%%%%%%%%%%%%%%%%%%%%%%%%%%%%%%%%%%%%%%%%
\medskip
\textbf{Acknowledgements.} {R.~H. gratefully acknowledges the
financial support of the Alexander von Humboldt Foundation and
thanks the Institute for Applied Mathematics of Bonn University
for the warm hospitality.}

%%%%%%%%%%%%%%%%%%%%%%%%%%%%%%%%%%%%%%%%%%%
\appendix

\section{Sobolev spaces $W^s_2(0,1)$ and all that}\label{sec:sobolev}

We recall here some facts about the Sobolev spaces $W^s_2(0,1)$
and Fourier coefficients of functions from these spaces. For
details, we refer the reader to~\cite[Ch.~1]{LM}.

By definition, the space $W^0_2(0,1)$ coincides with $L_2(0,1)$ and the
norm~$\|\cdot\|_0$ in $W^0_2(0,1)$ is just the $L_2(0,1)$-norm. The
Sobolev space $W^2_2(0,1)$ consists of all functions $f$ in $L_2(0,1)$,
whose distributional derivatives $f'$ and $f''$ also fall
into~$L_2(0,1)$. Being endowed with the norm
\[
    \|f\|_2 := \bigl(\|f\|_0^2 + \|f'\|_0^2 + \|f''\|_0^2\bigr)^{1/2},
\]
$W_2^2(0,1)$ becomes a Hilbert space.

Now we interpolate~\cite[CH.~1.2.1]{LM} between $W^2_2(0,1)$ and
$W^0_2(0,1)$ to get the intermediate spaces $W^s_2(0,1)$ with
norms $\|\cdot\|_s$ for $s\in(0,2)$; namely,
\[
    W^{2s}_2(0,1) :=
        [W^2_2(0,1),W^0_2(0,1)]_{1-s}.
\]
The norms $\|\cdot\|_s$ are nondecreasing with $s\in[0,2]$, i.e.,
if $s<t$ and $f\in W^t_2(0,1)$, then $\|f\|_s\le\|f\|_t$. Since by
construction the spaces $W_2^s(0,1)$ form an interpolation scale,
the general interpolation theorem~\cite[Theorem~1.5.1]{LM} implies
the following interpolation property for operators in these
spaces.

\begin{proposition}\label{pro:interp}
Assume that an operator $T$ acts boundedly in $W^s_2(0,1)$ and
$W^r_2(0,1)$, $s<r$. Then $T$ is a bounded operator in
$W^{ts+(1-t)r}_2(0,1)$ for every $t\in[0,1]$; moreover,
$\|T\|_{ts+(1-t)r}\le \|T\|_s^t \|T\|_r^{1-t}$.
\end{proposition}

Proposition~\ref{pro:interp} yields boundedness in every $W^s_2(0,1)$,
$s\in [0,2]$, of the operators $R$ and $V$ given by $Rf(x)=f(1-x)$ and
$Vf(x)=xf(x)$.

For an arbitrary $f\in L_2(0,1)$ and an arbitrary $\la\in\bC$, we put
\begin{equation}\label{eq:frtr}
    s_\la(f):= \int_0^1 f(x) \sin(\pi\la x)\,dx, \qquad
    c_\la(f):= \int_0^1 f(x) \cos(\pi\la x)\,dx.
\end{equation}
As usual, $\ast$ denotes the convolution operation on~$(0,1)$,
i.e.,
\[
    (f\ast g)(x) := \int_0^x f(x-t) g(t)\,dt.
\]
We shall also introduce the following shorthand notations:
\begin{align*}
    (f \oast g)(x) &:= \tfrac12 \bigl[ R\, (R f\ast g + f\ast Rg)
                + f\ast g + Rf \ast Rg\bigr],\\
    (f \hast g)(x) &:= \tfrac12 \bigl[ R\, (R f\ast g + f\ast Rg)
                - f\ast g - Rf \ast Rg\bigr],\\
    (f \tast g)(x) &:= \tfrac12 \bigl[ R\, (R f\ast g - f\ast Rg)
                + f\ast g - Rf \ast Rg\bigr]
\end{align*}
(where, as earlier, $R$ stands for the reflection operator, $(Rf)(x) =
f(1-x)$). The operations $\oast, \hast$, and $\tast$ play the same role for
the sine and cosine Fourier transform on $(0,1)$ as the usual convolution
for the Fourier transform on the whole line. Namely, these operations have
the following properties.

\begin{lemma}\label{lem:conv}
For arbitrary $f,g\in L_2(0,1)$ and $\la\in\bC$ the following equalities
hold:
\[
    c_\la(f)c_\la(g) = c_\la(f\oast g), \qquad
    s_\la(f)s_\la(g) = c_\la(f\hast g), \qquad
    s_\la(f)c_\la(g) = s_\la(f\tast g).
\]
\end{lemma}

\begin{proof}
We shall prove only the first equality since the other two can be treated
analogously. We have
\[
    2 c_\la(f)c_\la(g) = \int_0^1 \int_0^1 f(x) g(t)
            \{\cos[\pi \la(x-t)] + \cos[\pi \la(x+t)]\}\,dx dt,
\]
and simple calculations lead to
\begin{align*}
    \int_0^1\int_0^1 &f(x) g(t)\cos\pi \la(x-t)\, dx dt \\
        &= \int_0^1
        \Bigl( \int_0^{1-s} f(s+t) g(t)\,dt
                + \int_0^{1-s} f(t) g(s+t)\,dt
        \Bigr)\cos(\pi \la s)\,ds, \\
    \int_0^1\int_0^1 &f(x) g(t)\cos[\pi \la(x+t)]\, dx dt \\
        &= \int_0^1
        \Bigl( \int_0^{s} f(s-t) g(t)\,dt
                + \int_0^{s} f(1-t) g(1-s+t)\,dt
        \Bigr)\cos\pi \la s\,ds.
\end{align*}
Taking into account the relations
\[                          % 5.16
    \int_0^{1-s} f(s+t) g(t)\,dt = R(Rf\ast g)(s), \qquad
    \int_0^{s} f(1-t) g(1-s+t)\,dt = Rf \ast Rg,
\]
we get $c_\la(f) c_\la(g) = c_\la(f\oast g)$ as stated. The lemma is
proved.
\end{proof}

It is well known that convolution accumulates smoothness; the
precise meaning of this statement is as follows.

\begin{proposition}\label{pro:conv}
Assume that $s,t\in[0,1]$ and that $f\in W^s_2(0,1)$ and $g\in
W^t_2(0,1)$ are arbitrary. Then the function $h:= f\ast g$ belongs
to $W^{s+t}_2(0,1)$ and, moreover, there exists $C>0$ independent
of $f$ and $g$ such that $\|h\|_{s+t}\le C \|f\|_s\|g\|_t$.
\end{proposition}

Proof of this proposition is based on interpolation between the
extreme cases $s,t=0,1$, which are handled with directly.

Combining Proposition~\ref{pro:conv} with the fact that the
operator $R$ is bounded in the spaces~$W^s_2(0,1)$ for all
$s\in[0,1]$, we arrive at the following conclusion.

\begin{corollary}\label{cor:conv}
Assume that $s,t\in[0,1]$ and that $f\in W^s_2(0,1)$, $g\in W^t_2(0,1)$.
Then the functions $f\oast g, f\hast g$, and $f\tast g$ belong to
$W^{s+t}_2(0,1)$ and, moreover, there exists a number $C>0$ independent of
$f$ and $g$ such that
\[
    \max\bigl\{\|f\oast g\|_{s+t}, \|f\hast g\|_{s+t},
    \|f\tast g\|_{s+t} \bigr\} \le C \|f\|_s
    \|g\|_t.
\]
\end{corollary}

The following lemma is an analogue of the well-known Wiener lemma.

\begin{lemma}\label{lem:wiener}
Assume that $f\in W^s_2(0,1)$, where $s\in[0,1]$. If $1+ c_n(f)\ne0$ for
all $n\in\bN$, then there exists a function $g\in W^s_2(0,1)$ such that
\[
    \bigl(1 + c_n(f)\bigr)^{-1} = 1 + c_n(g), \qquad n\in\bN.
\]
\end{lemma}

\begin{proof}
We start the proof with some auxiliary constructions. Denote by $C$ the
operator that acts from $L_2(0,1)/\mathbb{C}$ into $\ell_2$ according
to the formula
\[
    Cf :=   \bigl( c_n(f)\bigr)_{n\in\bN}.
\]
This operator is isomorphic; we put
\[
    {\cW}^s:= \{ Cf \mid f\in W^s_2(0,1)\}
\]
and endow the linear space $\cW^s$ with the norm
\[
    \|x\|_{\cW^s}:= \|C^{-1}\,x\|_{s}.
\]
In view of Lemma~\ref{lem:conv} and Corollary~\ref{cor:conv}, the
elementwise multiplication $(xy)_n:=x_ny_n$ is a continuous
operation in~$\cW^s$. We adjoin to $\cW^s$ the unit element $e$
(with components $e_n$ equal to $1$) and denote the resulting
unital algebra by $\wh\cW^s$. By a well-known
result~\cite[Theorem~10.2]{Ru} one can introduce an equivalent
norm in $\wh\cW^s$ under which $\wh\cW^s$ becomes a commutative
Banach algebra.

Assume now that the assumptions of the lemma hold and denote by $x$ an
element of $\wh\cW^s$ with components $x_n:= 1+c_n(f)$. We shall prove
below that $x$ is invertible in $\wh\cW^s$; it then follows that $x^{-1} =
e + y$ for some $y\in \cW^s$ as required.

It is well known~\cite[Theorem~11.5]{Ru} that the element $x$ is
invertible in the unital Banach algebra~$\wh\cW^s$ if and only if
$x$ does not belong to any maximal ideal of $\wh\cW^s$. Assume, on
the contrary, that there exists a maximal ideal $\fm$ of
$\wh\cW^s$ containing $x$. Since $\wh\cW^s$ contains all finite
sequences and none of $x_n$ vanishes, $\fm$ also contains all
finite sequences. Finite sequences form a dense subset of $\cW^s$
because the set of all trigonometric polynomials in $\cos\pi nx$
is dense in $W^s_2(0,1)$. Recalling that maximal ideals are
closed, we conclude that $\cW^s\subset\fm$. Next we observe that
$\cW^s$ is a proper subset of $\fm$ (e.g., $x$ belongs to
$\fm\setminus\cW^s$) and that $\cW^s$ has codimension~$1$ in
$\wh\cW^s$. Henceforth $\fm = \wh\cW^s$, which contradicts our
assumption that $\fm$ is a maximal ideal of $\wh\cW^s$. As a
result, $x$ is not contained in any maximal ideal of $\wh\cW^s$
and thus is invertible in~$\wh\cW^s$. The lemma is proved.
\end{proof}

%%%%%%%%%%%%%%%%%%%%%%%%%%%%%%%%%%%%%%%%%%%

\end{document}